\documentclass[a4paper,10pt]{amsart}

\usepackage{amssymb}
\usepackage{bbm}
\usepackage{mathrsfs}	
\usepackage[T5,T1]{fontenc}

\normalsize

\newtheoremstyle{bracket}{1ex}{2ex}{\rm}{}{\bfseries}{}{0.8em}{\thmnumber{(#2)}}
\newtheoremstyle{thm}{1ex}{2ex}{\itshape}{}{\bfseries}{}{0.9em}{\thmnumber{(#2)}\thmname{ #1}\thmnote{ (#3)}}

\theoremstyle{bracket}
\newtheorem{no}{}

\theoremstyle{thm}
\newtheorem{prop}[no]{Proposition}
\newtheorem{theorem}[no]{Theorem}

\DeclareMathOperator{\spec}{Spec}
\DeclareMathOperator{\pic}{Pic}
\DeclareMathOperator{\ke}{Ker}

\newcommand{\dfgl}{\mathrel{\mathop:}=}

\newcommand{\ftnt}{\footnotetext{This is a slightly updated and corrected version of the author's contribution to the Proceedings of the 32nd Symposium and the 6th Japan-Vietnam Joint Seminar on Commutative Algebra, held in December 2010 in Hayama, Japan. The author was supported by the Swiss National Science Foundation.}}

\begin{document}

\title{On Toric Schemes\ftnt}
\author{Fred Rohrer}
\address{Institute of Mathematics, Vietnam Academy of Science and Technology, {\fontencoding{T5}\selectfont 18 Ho\`ang Qu\'\ocircumflex{}c Vi\d\ecircumflex{}t\\10307 H\`a N\d\ocircumflex{}i\\Vi\d\ecircumflex{}t Nam}}
\email{fredrohrer0@gmail.com}

\begin{abstract}
Studying toric varieties from a scheme-theoretical point of view leads to toric schemes, i.e.~``toric varieties over arbitrary base rings''. It is shown how the base ring affects the geometry of a toric scheme. Moreover, generalisations of results by Cox and Musta\c{t}\v{a} allow to describe quasicoherent sheaves on toric schemes in terms of graded modules. Finally, a toric version of the Serre-Grothendieck correspondence relates cohomology of quasicoherent sheaves on toric schemes to local cohomology of graded modules.
\end{abstract}

\maketitle

\setcounter{section}{-1}
\setcounter{no}{-1}


\section{From toric varieties to toric schemes}

During the last forty years a huge amount of work on toric varieties was and still is published. Their theory was generalised in several directions, and this often lead to a better understanding of classical toric varieties. However, the generalisation that seems to be the most natural and the most important -- the \textit{study of toric varieties from a scheme-theoretical point of view} -- was never actually carried out. It is clear that to do this one has to be able to make arbitrary base changes. Hence, instead of considering toric varieties over an algebraically closed field (or, as often done, over the field of complex numbers), one needs to study \textit{toric schemes,} that is ``toric varieties over arbitrary base rings''. Special cases of this generalisation were mentioned briefly in \cite[\textsection 4]{dem} (for regular fans and mainly over the ring of integers) and \cite[IV.3]{kkms} (over discrete valuation rings). But unfortunately later authors seemed to ignore this, and hence the knowledge of toric schemes is very small compared to the one of toric varieties.

Besides yielding a better understanding of the geometry of toric varieties, there are concrete applications of the above generalisation, as the following remark shows.

\begin{no}
Let $X$ be the toric variety over an algebraically closed field $K$ associated with a fan $\Sigma$. A fundamental question in algebraic geometry is then if the Hilbert functor ${\rm Hilb}_{X/K}$ of $X$ over $K$ is representable, i.e. if the Hilbert scheme of $X$ exists (cf.~\cite{fga}). If $X$ is projective, then this is indeed the case and follows from Grothendieck's more general result \cite[Th\'eor\`eme 3.1]{fga}. However, toric varieties are not necessarily projective, and in general it is not known whether their Hilbert schemes exist. Studying ${\rm Hilb}_{X/K}$ amounts to studying quasicoherent sheaves on the base change $X\otimes_KR$ for every $K$-algebra $R$, and it turns out that $X\otimes_KR$ is the same as the toric scheme over $R$ associated with $\Sigma$. Hence, in order to study Hilbert functors of toric varieties \textit{it is necessary to study toric schemes over more general bases than just over algebraically closed fields.}
\end{no}

The development of a theory of toric schemes was begun in the PhD Thesis \cite{diss}, and its contents were refined and extended in \cite{cof}, \cite{ts1} and \cite{qcs}. Here we give an overview of the most important results and refer the reader to the aforementioned sources for a more extensive treatment including proofs.


\section{The geometry of toric schemes}

We start by briefly describing the construction of toric schemes from fans. In \cite{ts1}, toric schemes are obtained as a special case of the more general construction of schemes from so-called openly immersive projective systems of monoids (also yielding the Cox schemes introduced below).

\smallskip

\noindent\textit{$\bullet$\quad From now on let $V$ be an $\mathbbm{R}$-vector space of finite dimension $n$, let $N$ be a $\mathbbm{Z}$-structure on $V$ (i.e.~a subgroup of rank $n$ of the additive group underlying $V$ with $\langle N\rangle_{\mathbbm{R}}=V$), and let $M\dfgl N^*$ denote the dual of $N$ which is a $\mathbbm{Z}$-structure on the dual $V^*$ of $V$.}

\smallskip

An \textit{$N$-polycone (in $V$)} is the set of $\mathbbm{R}$-linear combinations with coefficients in $\mathbbm{R}_{\geq 0}$ of a finite subset of $N$, and an $N$-polycone is called \textit{sharp} if it does not contain a line. If $\sigma$ is an $N$-polycone then a \textit{face of $\sigma$} is a set of the form $\sigma\cap\ke(u)$ for some $u\in\sigma^{\vee}\cap M$ (where $E^{\vee}\dfgl\{v\in V^*\mid v(E)\subseteq\mathbbm{R}_{\geq 0}\}$ for a subset $E\subseteq V$). The set of faces of a (sharp) $N$-polycone is a finite set of (sharp) $N$-polycones. An \textit{$N$-fan (in $V$)} is a finite set $\Sigma$ of sharp $N$-polycones that is closed under taking faces and such that the intersection of two cones in $\Sigma$ is a common face of both of them. By means of the relation ``$\tau$ is a face of $\sigma$'', denoted by $\tau\preccurlyeq\sigma$, we consider an $N$-fan as an ordered set.

\smallskip

\noindent\textit{$\bullet$\quad From now on let $\Sigma$ be an $N$-fan in $V$ and let $R$ be a ring\/\footnote{By a ring, group or monoid we always mean a commutative ring, group or monoid, respectively, and by an algebra we always mean a commutative, unital and associative algebra.}.}

\smallskip

If $\sigma\in\Sigma$ then $\sigma^{\vee}\cap M$ is a torsionfree, cancellable, finitely generated submonoid of $M$, and if moreover $\tau\preccurlyeq\sigma$ then $\sigma^{\vee}\cap M$ is a submonoid of $\tau^{\vee}\cap M$. Taking spectra of algebras of monoids over $R$ and setting $X_{\sigma}(R)\dfgl\spec(R[\sigma^{\vee}\cap M])$ for $\sigma\in\Sigma$ we get an inductive system $(X_{\sigma}(R))_{\sigma\in\Sigma}$ of $R$-schemes over $\Sigma$. Its inductive limit exists and is an $R$-scheme denoted by $X_{\Sigma}(R)\rightarrow\spec(R)$ and called \textit{the toric scheme over $R$ associated with $\Sigma$ (and $N$).} It can be understood as obtained by glueing $(X_{\sigma}(R))_{\sigma\in\Sigma}$ along $(X_{\sigma\cap\tau}(R))_{(\sigma,\tau)\in\Sigma^2}$.

The above construction of toric schemes gives rise to a contravariant functor $X_{\Sigma}$ from the category of rings to the category of schemes together with a morphism $X_{\Sigma}\rightarrow\spec$. Moreover, the functor $X_{\Sigma}$ is compatible with base change in the following sense.

\begin{prop}\mbox{\rm(\cite[1.6]{ts1})}\label{basechange}
There is a canonical isomorphism $$X_{\Sigma}(\bullet)\cong X_{\Sigma}(R)\otimes_R\bullet$$ of contravariant functors from the category of $R$-algebras to the category of $R$-schemes.
\end{prop}

In particular, if $\mathfrak{a}\subseteq R$ is an ideal then $X_{\Sigma}(R/\mathfrak{a})$ is canonically identified with a closed subscheme of $X_{\Sigma}(R)$.

\smallskip

The first important question is now of course how the base ring affects the geometry of a toric scheme. It turns out that some basic properties hold for all toric schemes, making them a class of ``nice schemes''. More precisely, on use of the above base change property we get the following result.

\begin{prop}\mbox{\rm(\cite[3.4]{ts1})}
The $R$-scheme $X_{\Sigma}(R)\rightarrow\spec(R)$ is separated, quasicompact, flat, and of finite presentation; it is faithfully flat if and only if\/ $\Sigma\neq\emptyset$ or $R=0$.
\end{prop}

In contrast, a lot of other basic properties are respected and reflected by $X_{\Sigma}$. The following statements are proved by reducing to the affine case, i.e.~$X_{\sigma}$, and then applying corresponding results about algebras of monoids (see e.g.~\cite{gil}).

\begin{prop}\mbox{\rm(\cite[3.4]{ts1})}
a) The scheme $X_{\Sigma}(R)$ is reduced, connected, or normal if and only if $R$ is so or\/ $\Sigma=\emptyset$; it is irreducible, or integral if and only if $R$ is so and\/ $\Sigma\neq\emptyset$.

b) If\/ $\Sigma\neq\emptyset$ then there is a bijection\/ $\mathfrak{p}\mapsto X_{\Sigma}(R/\mathfrak{p})$ from the set of minimal prime ideals of $R$ to the set of irreducible components of $X_{\Sigma}(R)$.

c) The scheme $X_{\Sigma}(R)$ is Noetherian if and only if $R$ is so or\/ $\Sigma=\emptyset$; it is Artinian if and only if $R$ is so and $n=0$, or $R=0$, or\/ $\Sigma=\emptyset$.

d) If\/ $\Sigma\neq\emptyset$ then $$\dim(R)+n\leq\dim(X_{\Sigma}(R))\leq(n+1)\dim(R)+n;$$ if $R$ is moreover Noetherian then $\dim(R)+n=\dim(X_{\Sigma}(R))$.

e) If $R$ is Noetherian, then $X_{\Sigma}(R)$ is equidimensional if and only if $R$ is so or\/ $\Sigma=\emptyset$.
\end{prop}

The above shows in particular that on general toric schemes \textit{no satisfying theory of Weil divisors is available.} Since a lot of results about toric varieties were proved by heavy use of Weil divisor techniques (see e.g.~\cite{cox}, \cite{ful}), one has to come up with new proofs in order to generalise these results to toric schemes.

\smallskip

Finally, as an example of a property depending on the fan but not on the base ring we consider properness. Its characterisation needs the notion of a \textit{complete} $N$-fan $\Sigma$, i.e.~an $N$-fan $\Sigma$ with $\bigcup\Sigma=V$. This result is well-known for toric varieties (see e.g.~\cite[2.4]{ful}), and proved on use of torus operations for toric schemes associated with regular fans in \cite[\textsection 4 Proposition 4]{dem}. Our proof for arbitrary fans avoids speaking of torus operations and relies only on the valuative criterion for properness and on properties of projections of fans proved in \cite{cof}.

\begin{prop}
The $R$-scheme $X_{\Sigma}(R)\rightarrow\spec(R)$ is proper if and only if\/ $\Sigma$ is complete, or\/ $\Sigma=\emptyset$, or $R=0$.
\end{prop}


\section{Sheaves on toric schemes}

Generalising work by Cox \cite{cox} and Musta\c{t}\v{a} \cite{mus1} we introduce a notion of Cox ring (not to be confused with the one introduced in \cite{hukeel}) and describe quasicoherent modules on toric schemes in terms of graded modules over these rings. In order to do so we need to define some objects encoding the combinatorics of the fan $\Sigma$.

\smallskip

Let $\Sigma_1$ denote the set of $1$-dimensional cones in $\Sigma$. Every $\rho\in\Sigma_1$ has a unique minimal $N$-generator (i.e.~an $x\in N$ with $\rho=\mathbbm{R}_{\geq 0}x$ such that $rx\notin N$ for every $r\in\ ]0,1[$), denoted by $\rho_N$. There is an exact sequence of groups $$M\overset{c}\longrightarrow\mathbbm{Z}^{\Sigma_1}\overset{a}\longrightarrow A\longrightarrow 0,$$ where $c(u)\dfgl(u(\rho_N))_{\rho\in\Sigma_1}$ for $u\in M$ and where $a$ is defined as the cokernel of $c$. Note that $c$ is a monomorphism if and only if $\Sigma$ if \textit{full,} i.e.~$\langle\bigcup\Sigma\rangle_{\mathbbm{R}}=V$. We denote by $(\delta_{\rho})_{\rho\in\Sigma_1}$ the canonical basis of $\mathbbm{Z}^{\Sigma_1}$ and we set $\alpha_{\rho}\dfgl a(\delta_{\rho})$ for $\rho\in\Sigma_1$.

Now, we denote by $S$ the polynomial algebra $R[(Z_{\rho})_{\rho\in\Sigma_1}]$ in indeterminates $(Z_{\rho})_{\rho\in\Sigma_1}$ over $R$, furnished with the $A$-graduation induced by $a$, i.e.~such that $\deg(Z_{\rho})=\alpha_{\rho}$ for $\rho\in\Sigma_1$. For $\sigma\in\Sigma$ we set $\widehat{Z}_{\sigma}\dfgl\prod_{\rho\in\Sigma_1\setminus\sigma_1}Z_{\rho}\in S$ (where $\sigma_1$ denotes the set of $1$-dimensional faces of $\sigma$). Finally we define a graded ideal $I\dfgl\langle\widehat{Z}_{\sigma}\mid\sigma\in\Sigma\rangle_S$.

\smallskip

\noindent\textit{$\bullet$\quad From now on let $B\subseteq A$ be a subgroup.}

\smallskip

The $B$-graded $R$-algebra $S_B\dfgl\bigoplus_{\alpha\in B}S_{\alpha}$ obtained from $S$ by degree restriction to $B$ is called \textit{the $B$-restricted Cox ring over $R$ associated with $\Sigma_1$ (and $N$),} and its graded ideal $I_B\dfgl I\cap S_B$ is called \textit{the $B$-restricted irrelevant ideal over $R$ associated with $\Sigma$ (and $N$).} One can show that $I_B$ is generated by finitely many monomials.

To proceed we need to ``invert the monomials $\widehat{Z}_{\sigma}$ in the Cox ring'', and hence we have to assure that some power of these monomials lies in $S_B$. This amounts to supposing that $B$ is \textit{big,} i.e.~it has finite index in $A$.

\smallskip

\noindent\textit{$\bullet$\quad From now on suppose that $B$ is big, so that there exists $m\in\mathbbm{N}_0$ with $\widehat{Z}_{\sigma}^m\in S_B$ for every $\sigma\in\Sigma$.}

\smallskip

For $\sigma\in\Sigma$ the $B$-graded ring of fractions $(S_B)_{\widehat{Z}_{\sigma}^m}$ is independent of the choice of $m$. Its component of degree $0$ is independent of the choice of $B$ and is denoted by $S_{(\sigma)}$. Moreover, for $\tau\preccurlyeq\sigma$ there is a canonical morphism of rings $S_{(\sigma)}\rightarrow S_{(\tau)}$ which is independent of $m$ and $B$. Taking spectra and setting $Y_{(\sigma)}(R)\dfgl\spec(S_{(\sigma)})$ for $\sigma\in\Sigma$ we obtain an inductive system $(Y_{\sigma}(R))_{\sigma\in\Sigma}$ of $R$-schemes over $\Sigma$. Its inductive limit exists and is an $R$-scheme denoted by $Y_{\Sigma}(R)\rightarrow\spec(R)$ and called \textit{the Cox scheme over $R$ associated with $\Sigma$ (and $N$).} It can be understood as obtained by glueing $(Y_{\sigma}(R))_{\sigma\in\Sigma}$ along $(Y_{\sigma\cap\tau}(R))_{(\sigma,\tau)\in\Sigma^2}$.

The above construction of Cox schemes gives rise to a contravariant functor $Y_{\Sigma}$ from the category of rings to the category of schemes together with a morphism $Y_{\Sigma}\rightarrow\spec$, and $Y_{\Sigma}$ is compatible with base change in the sense of (\ref{basechange}).

\smallskip

Cox schemes are closely related to toric schemes as follows. The morphism of groups $c\colon M\rightarrow\mathbbm{Z}^{\Sigma_1}$ induces morphisms of rings $R[\sigma^{\vee}\cap M]\rightarrow S_{(\sigma)}$ for $\sigma\in\Sigma$, and these induce a canonical morphism of contravariant functors $\gamma\colon Y_{\Sigma}\rightarrow X_{\Sigma}$. Then, we have the following result.

\begin{prop}\mbox{\rm(\cite[3.3.3]{qcs})}
The canonical morphism of contravariant functors $\gamma\colon Y_{\Sigma}\rightarrow X_{\Sigma}$ is an isomorphism if and only if\/ $\Sigma$ is full.
\end{prop}

Using the (non-canonical) procedure to consider a toric scheme associated with a non-full fan as a toric scheme associated with a full fan (\cite[3.3]{ts1}) it is sufficient to study from now on Cox schemes instead of toric schemes. (Note that this reduction demands a base change and is in general \textit{not available for toric varieties.})

\medskip

Now we are ready to explain how $B$-graded $S_B$-modules give rise to quasicoherent sheaves on $Y_{\Sigma}(R)$. We denote by ${\sf GrMod}^B(S_B)$ and ${\sf QCMod}(\mathscr{O}_{Y_{\Sigma}(R)})$ the categories of $B$-graded $S_B$-modules and of quasicoherent $\mathscr{O}_{Y_{\Sigma}(R)}$-modules. Moreover, for a $B$-graded $S_B$-module $F$ we denote by $F_{(\sigma)}$ the component of degree $0$ of the $B$-graded module of fractions $F_{\widehat{Z}_{\sigma}^m}=F\otimes_{S_B}(S_B)_{\widehat{Z}_{\sigma}^m}$, and for an $S_{(\sigma)}$-module $G$ we denote by $\widetilde{G}$ the $\mathscr{O}_{Y_{\sigma}(R)}$-module associated with $G$.

\begin{prop}\mbox{\rm(\cite[4.1.1]{qcs})}
There exists a unique functor $$\mathscr{S}_B\colon{\sf GrMod}^B(S_B)\rightarrow{\sf QCMod}(\mathscr{O}_{Y_{\Sigma}(R)})$$ with $\mathscr{S}_B(F)\!\upharpoonright_{Y_{\sigma}(R)}=\widetilde{F_{(\sigma)}}$ for every $\sigma\in\Sigma$ and every $B$-graded $S_B$-module $F$.
\end{prop}

Since $\mathscr{S}_B$ coincides locally with the canonical equivalence between modules and quasicoherent sheaves on affine schemes it is exact and commutes with inductive limits. Furthermore, denoting by $\bullet(\alpha)$ the functor of shifting degrees by $\alpha$, we can construct a right quasiinverse $$\Gamma_*^B(\bullet)\dfgl\bigoplus_{\alpha\in B}\Gamma\bigl(Y_{\Sigma}(R),\bigl(\bullet\otimes_{\mathscr{O}_{Y_{\Sigma}(R)}}\mathscr{S}_B(S_B(\alpha))\bigr)\bigr)$$ for $\mathscr{S}_B$, called \textit{the first total functor of sections associated with $\Sigma$ and $B$ over $R$.} Thus, we get the following generalisation of \cite[Theorem 1.1]{mus1}, itself a generalisation of \cite[Theorem 3.2]{cox}.

\begin{theorem}\label{surj}\mbox{\rm(\cite[4.4.3]{qcs})}
The functor $\mathscr{S}_B\colon{\sf GrMod}^B(S_B)\rightarrow{\sf QCMod}(\mathscr{O}_{Y_{\Sigma}(R)})$ is essentially surjective.
\end{theorem}

Next, we restrict our attention to ideals. A graded ideal $\mathfrak{a}\subseteq S_B$ is called \textit{$I_B$-saturated} if $\mathfrak{a}=\bigcup_{k\in\mathbbm{N}_0}(\mathfrak{a}:_{S_B}I_B^k)$. Let $\mathbbm{J}^{{\rm sat}}_B$ and $\widetilde{\mathbbm{J}}$ denote the sets of $I_B$-saturated graded ideals of $S_B$ and of quasicoherent ideals of $\mathscr{O}_{Y_{\Sigma}(R)}$, respectively. Then, $\mathscr{S}_B$ induces by exactness a map $\Xi_B\colon\mathbbm{J}^{{\rm sat}}_B\rightarrow\widetilde{\mathbbm{J}}$. The next result treats the question whether this map is surjective or injective. To get injectivity, besides being big the subgroup $B$ must not be ``too big''. More precisely, $B$ is called \textit{small (with respect to $\Sigma$)} if it is contained in $\bigcap_{\sigma\in\Sigma}\langle\{\alpha_{\rho}\mid\rho\in\Sigma_1\setminus\sigma_1\}\rangle_{\mathbbm{Z}}$.

\begin{theorem}\label{bij}\mbox{\rm(\cite[4.4.9]{qcs})}
The map\/ $\Xi_B\colon\mathbbm{J}_B^{{\rm sat}}\rightarrow\widetilde{\mathbbm{J}}$ is surjective, and if $B$ is small then it is bijective.
\end{theorem}

An example of a subgroup that is big and small (and moreover well understood) is given in the following remark (cf.~\cite[V.5]{ewald}).

Consider a family $(U_{\sigma})_{\sigma\in\Sigma}$ of subsets of $V^*$ such that for every $\sigma\in\Sigma$ there exists a (not necessarily unique) $m_{\sigma}\in M$ with $U_{\sigma}=m_{\sigma}+\sigma^{\vee}$. Such a family is called a \textit{virtual polytope over\/ $\Sigma$} if $\tau\subseteq\ke(m_{\sigma}-m_{\tau})$ for all $\sigma,\tau\in\Sigma$ with $\tau\preccurlyeq\sigma$, and this condition is independent of the choice of the family $(m_{\sigma})_{\sigma\in\Sigma}$. There is a canonical structure of group on the set of virtual polytopes over $\Sigma$, and the set of virtual polytopes of the form $(m+\sigma^{\vee})_{\sigma\in\Sigma}$ is a subgroup. The corresponding quotient group is denoted by $\pic(\Sigma)$ and called \textit{the Picard group of\/ $\Sigma$.} It can be considered as the group of virtual polytopes over $\Sigma$ modulo $M$-rational translations.

The map $$(m_{\sigma}+\sigma^{\vee})_{\sigma\in\Sigma}\mapsto(m_{\rho}(\rho_N))_{\rho\in\Sigma_1}$$ yields a monomorphism from the group of virtual polytopes over $\Sigma$ to $\mathbbm{Z}^{\Sigma_1}$, and this induces a monomorphism $\pic(\Sigma)\rightarrowtail A$ by means of which we consider $\pic(\Sigma)$ as a subgroup of $A$. Then, $\pic(\Sigma)$ is small, and if $\Sigma$ is simplicial then $\pic(\Sigma)$ is big. Hence, it provides an example of a subgroup of $A$ to which (\ref{bij}) can be applied.

Finally, since $\pic(\Sigma)\cong\pic(X_{\Sigma}(\mathbbm{C}))$ by \cite[Theorem VII.2.15]{ewald} we get back \cite[Corollary 3.9]{cox} as a special case.


\section{Cohomology on toric schemes}

Our results about quasicoherent sheaves in the last section reveals that toric schemes (or more precisely, Cox schemes) are very similar to projective schemes. Hence, we ask if there is a toric version of the Serre-Grothendieck correspondence (cf.~\cite[20.4.4]{bs}), relating cohomology of quasicoherent sheaves on a Cox scheme to graded local cohomology of $B$-graded $S_B$-modules with respect to the irrelevant ideal $I_B$. This is indeed the case.

\smallskip

First, we have to explain what we mean by graded local cohomology. We denote by $${}^B\Gamma_{I_B}\colon{\sf GrMod}^B(S_B)\rightarrow{\sf GrMod}^B(S_B)$$ the $B$-graded $I_B$-torsion functor. Its right derived cohomological functor is denoted by $({}^BH^i_{I_B})_{i\in\mathbbm{Z}}$ and called \textit{$B$-graded local cohomology with respect to $I_B$.} The reason for this clumsy notation is that the ungraded module underlying a graded local cohomology module of a graded module $F$ might not be the same as the local cohomology module of the ungraded module underlying $F$. (A sufficient condition for this to hold is coherence of the graded ring $S_B$.)

Next, we introduce a variant of sheaf cohomology that is useful for our purpose. We define a functor $$\Gamma_{**}^B(\bullet)\colon{\sf GrMod}^B(S_B)\rightarrow{\sf GrMod}^B(S_B),$$ called \textit{the second total functor of sections associated with $\Sigma$ and $B$ over $R$,} by setting $$\Gamma_{**}^B(\bullet)\dfgl\bigoplus_{\alpha\in B}\Gamma(Y_{\Sigma}(R),\mathscr{S}_B(\bullet(\alpha))).$$ Note that despite its name it is defined on the category ${\sf GrMod}^B(S_B)$. However, by (\ref{surj}) this is merely a technical point. The reason for two (in general different) total functors of sections is that the canonical morphism $$\mathscr{S}_B(\bullet)\otimes_{\mathscr{O}_{Y_{\Sigma}(R)}}\mathscr{S}_B(S_B(\alpha))\rightarrow\mathscr{S}_B(\bullet(\alpha))$$ is not necessarily an isomorphism. The right derived cohomological functor of $\Gamma_{**}^B(\bullet)$ is denoted by $(H^i_{**,B})_{i\in\mathbbm{Z}}$ and contains the usual sheaf cohomology as a direct summand.

To go on we need a certain behaviour of injectives in the category ${\sf GrMod}^B(S_B)$. Namely, the $B$-graded ring $S_B$ is said to have \textit{the ITR-property with respect to $I_B$} if every $B$-graded $I_B$-torsion $S_B$-module has an injective resolution whose components are $B$-graded $I_B$-torsion $S_B$-modules. This is fulfilled for example if $S_B$ is Noetherian (as a graded ring), and in particular if $R$ is Noetherian. Using this notion and imitating the corresponding proof in the projective case we arrive at the Toric Serre-Grothendieck Correspondence.

\begin{theorem}\mbox{\rm(\cite[4.5.4]{qcs})}
If $S_B$ has the ITR-property with respect to $I_B$, then there exist an exact sequence of functors $$0\longrightarrow{}^B\Gamma_{I_B}\longrightarrow{\rm Id}_{{\sf GrMod}^B(S_B)}\longrightarrow\Gamma_{**}^B\overset{\zeta_B}\longrightarrow{}^BH^1_{I_B}\longrightarrow 0$$ and a unique morphism of $\delta$-functors $$(\zeta^i_B)_{i\in\mathbbm{Z}}\colon\bigl(H^i_{**,B}\bigr)_{i\in\mathbbm{Z}}\longrightarrow\bigl({}^BH^{i+1}_{I_B}\bigr)_{i\in\mathbbm{Z}}$$ with $\zeta^0_B=\zeta_B$, and $\zeta^i_B$ is an isomorphism for every $i\in\mathbbm{N}$.
\end{theorem}

As an application we can prove a toric version of Serre's Finiteness Theorem.

\begin{prop}\mbox{\rm(\cite[4.5.5]{qcs})}
Let $F$ be a finitely generated $B$-graded $S_B$-module, and suppose that $\Sigma$ is complete and that $R$ is Noetherian. Then, the $R$-modules $H^i_{**,B}(F)_{\alpha}$ and ${}^BH^i_{I_B}(F)_{\alpha}$ are finitely generated for every $i\in\mathbbm{Z}$ and every $\alpha\in B$.
\end{prop}

Considering the fibres of a toric scheme, this allows us to define and investigate Hilbert functions of toric schemes, a task we would like to address in future research. Note that the above hypothesis of a complete fan $\Sigma$ can be achieved by the Completion Theorem (\cite{cof}).


\end{document}